\documentclass{amsart}
    \usepackage{array}
    \usepackage{amsmath}
    \usepackage{amsfonts,amssymb}
\usepackage{amsthm}

\newcolumntype{C}{>{$}c<{$}}

\newtheorem{example}{Example}[section]
\newtheorem{theorem}{Theorem}[section]
\newtheorem{definition}{Definition}[section]
\newtheorem{proposition}{Proposition}[section]
\newtheorem{remark}{Remark}[section]

\newtheorem{corollary}{Corollary}[section]

\def\gc{\mathfrak g_{{}_{\mathbf C}}}
\def\kc{\mathfrak k_{{}_{\mathbf C}}}
\def\pc{\mathfrak p_{{}_{\mathbf C}}}

\def\qc{{\mathfrak q}_{{}_{\mathbf C}}}
\def\lc{\mathfrak l_{{}_{\mathbf C}}}

\def\nk{{\mathfrak n}_k}

\def\lk{\lc\cap\kc}

\def\Kc{K_{{}_{\mathbf C}}}
\def\Gc{G_{{}_{\mathbf C}}}

\def\Uc{U_{{}_{\mathbf C}}}
\def\Lc{L_{{}_{\mathbf C}}}

\def\Qc{{\mathbf Q}_{{}_{\mathbf C}}}
\def\tc{\mathfrak t_{{}_{\mathbf C}}}

\def\uc{\mathfrak u_{{}_{\mathbf C}}}

\def\g{\mathfrak g}
\def\k{\mathfrak k}

\def\t{\mathfrak t}

\begin{document}

\title{Small spherical nilpotent orbits and $K$-types of Harish Chandra modules}
\author{\text{Donald R. King}\\ $\textit{Department of Mathematics}$\\ $\textit{Northeastern University}$\\ $\textit{Boston,
Massachusetts 02115}$\\ $\text{E-mail: d.king@neu.edu}$}
\date{}
%\subjclass{Primary 22E46; Secondary 14R20, 53D20.}

\begin{abstract} Let $G$ be a connected linear semisimple Lie group with Lie algebra $\mathfrak g$ and maximal compact subgroup $K$. Let ${\Kc}~\rightarrow~{Aut(\pc)}$ be the
complexified isotropy representation at the identity coset of the
corresponding symmetric space. Suppose that $\mathcal O$ is a
nilpotent $\Kc$-orbit in $\pc$, and $\overline{\mathcal O}$ is its
Zariski closure in $\pc$. We study the $K$-type decomposition of
the ring of regular functions on $\overline{\mathcal O}$ when
$\mathcal O$ is spherical and ``small''. We also show that this
decomposition gives the asymptotic directions of $K$-types in any
irreducible $(\gc,\ K)$-module whose associated variety is
$\overline{\mathcal O}$.

\thanks{2000 {\it Mathematics Subject Classificationn}. Primary 22E46; Secondary 14R20, 53D20.}

\end{abstract}
\maketitle
\markboth{D. R. King}{Ring of regular functions}

\section{A desingularisation of the closure of a nilpotent $\Kc$-orbit in $\pc$}
    Let $G$ be a connected linear semisimple Lie group with Lie
algebra $\mathfrak g$ and maximal compact subgroup $K$. Let
${\Kc}~\rightarrow~{Aut(\pc)}$ be the complexified isotropy
representation at the identity coset of the corresponding
symmetric space. Suppose that $\mathcal O$ is a nilpotent
$\Kc$-orbit in $\pc$, and $\overline{\mathcal O}$ is its Zariski
closure in $\pc$. ${\overline{\mathcal O}}^n$ denotes the
normalization of ${\overline{\mathcal O}}$. If $W$ is an algebraic
variety, $R[W]$ denotes the ring of regular functions on
    $W$.\par

\medskip
    We recall some facts related to vector bundles over a
homogeneous space ${M~{\slash~P}}$ for a complex reductive group
$M$ and parabolic subgroup $P$. Let $V$ be a finite dimensional
$P$-module. Then, $M\times_{Q} V = {M\times V}\slash \sim$, where
the equivalence relation  ``$\sim$'' is defined by: $(gy,\ v)\sim
(g,\ y\cdot v)$, for $g\in M$, $y\in P$ and $v\in V$. The
equivalence relation ``$\sim$'' can be seen as arising from a
right action of $P$ on $M\times V$ defined as follows: $(g,\
v)\cdot y = (gy,\ y^{-1}\cdot v)$. Then it is easy to check that
$(g,\ v)\cdot (y_1y_2) = [(g,\ v)\cdot y_1]\cdot y_2$. Also, note
that $(g,\ v)\sim (g,\ v)\cdot y$, for all $g\in M,\ v\in V,\ y\in
P$. The equivalence class of $(g,\ v)$ under $\sim$ is denoted
$[g,\ v]$.\par
     A key fact is that:
\begin{equation*}
R[M\times_{P} V] = R[M\times V]^P,
\end{equation*}
where the $P$ action on $R[M\times V]$ is defined from the $P$
action on $M\times V$: if $k\in R[M\times V]$, and $y\in P$, then
$(k\cdot y)(g,\ v)= k((g,\ v)\cdot y^{-1})=k(gy^{-1},\ y\cdot v)$.
The equality of the rings above is implemented as follows. If
$f\in R[M\times_{P} V]$, then define $\widetilde{f}$ on $M\times
V$, by $\widetilde{f}(g,\ v) = f([g,\ v])$. Since $f([g,\ v]) =
f([gyy^{-1},\ v]) = f([gy,\ y^{-1}\cdot v])$, $\widetilde{f}\in
R[M\times V]^P$. On the other hand if $h\in R[M\times V]^P$, then
$\hat{h}([g,\ v]):=h(g,\ v)$ is well defined and $\hat{h}\in
R[M\times_{P} V]$.\par Furthermore,
\begin{equation*}
\begin{split}
&R[M\times V]^P=(R[M]\otimes_{\mathbf C}R[V])^P=
(R[M]\otimes_{\mathbf C}S[V^*])^P\\ &=\sum_{n\ge
0}(R[M]\otimes_{\mathbf C}S^n[V^*])^P=\sum_{n\ge 0}H^0(M\slash P,\
S^n(V^*)),
\end{split}
\end{equation*}
where $S(V^*)$ is the symmetric algebra on $V^*$, the dual of $V$.
Here, the action of $P$ on $R[M]\otimes_{\mathbf C}R[V]$ is:
$[(f_1(\cdot)\otimes f_2(\cdot))\cdot y](g,\ v) =
f_1(gy^{-1})f_2(y\cdot v)$

  \par
\bigskip
    Let $\{x,\ e, \ f\}$ be a normal triple, in the sense of Kostant and Rallis \cite{colmcg}, corresponding to the nilpotent $\Kc$-orbit $\mathcal O$ in $\pc$.
\begin{definition} (The eigenspaces of $ad(x)$ on $\gc$) Let $\gc(x;j)$
=  the $j$-eigenspace of $x$. Likewise define $\pc(x;j)$ and
$\kc(x;j)$. $\mathcal O$ will be said to have height $j$, if $\gc(x;j)\ne 0$ but $\gc(x;j')=0$ for all $j'>j$.
\end{definition}
\begin{definition}\label{desingdef} (Construction of a desingularization of $\overline{\mathcal O}$ .)  Set $V = V(e) = \sum^{}_{j\ge 2}\gc(x;j)$. $\qc$ =
$\sum\limits_{j\ge 0}\gc(x;j)$, $\uc= \sum\limits_{j>0}\gc(x;j)$,
and $\lc = \gc(x;0)$. Let $\Qc, \ \Lc$ , and $\Uc$ be the connected
subgroups of $\Gc$ with Lie algebras $\qc,\ \lc,\ {\rm and} \ \uc$
respectively. It is well known that the morphism $\pi
:\Gc\times_{\Qc} V \rightarrow \overline{\Gc\cdot{e}}$ , defined by
$\pi ([g,v)] = g\cdot v$ is a desingularization of
$\overline{\Gc\cdot{e}}$ . By similar arguments, if $\widetilde{V} =
\sum\limits_{j\ge 2}\pc(x;j)= V\cap \pc$, then the (restriction)
mapping
\begin{equation}\label{desingmathcalo}
\pi :\Kc\times _{\Qc\cap \Kc}\widetilde{V} \rightarrow
\overline{\Kc\cdot{e}}
\end{equation}
is a desingularization (resolution of singularities) of
$\overline{\mathcal O}$.
\end{definition}
\begin{remark}\label{gyyaapp}
Since $ad(e): \qc \rightarrow  V$ is surjective, $ad(e): \qc\cap\kc
\rightarrow  V\cap\pc $ is surjective. That is, $[\qc\cap\kc
,\ e]=V\cap\pc $. This enables us to conclude
that the orbit $Q\cap K_{{\Bbb C}}\cdot e$ is dense in $[\widetilde{V}$.
The remainder of the argument is essentially the same as the one
establishing that $\pi :\Gc\times _{\Qc} V \rightarrow
\overline{\Gc\cdot{e}}$  is a desingularization of
$\overline{\Gc\cdot{e}}$.)
\end{remark}
\begin{proposition}\label{rkctimes} With notation as above, $R[\Kc\times _{\Qc\cap \Kc}\widetilde{V}]=R[{\overline{\mathcal O}}^n]$
\end{proposition}

\section{Normality of closures of small spherical orbits}
\begin{definition} The nilpotent orbit $\mathcal O=\Kc\cdot e$ is said to be
spherical if some Borel subgroup of $\Kc$ has a dense orbit in
$\mathcal O$.
\end{definition}
    A classification of spherical nilpotent orbits for $\g$ real and simple may be found in
    \cite{king}.
\begin{definition}\label{smallspher} A spherical orbit $\mathcal O$ is said to be small if $\pc(x;i)=0$ for all $i>2$.
\end{definition}
\begin{remark} If an orbit $\mathcal O$ satisfies $\pc(x;i)=0$ for $i>2$, we have the exact sequence:
\begin{equation*}
0\rightarrow {\kc}^{\{x,\ e,\ f\}}\rightarrow \kc^x\rightarrow
\pc(x;2)\rightarrow 0,
\end{equation*}
 where the $ad(e):\kc^x\rightarrow \pc(x;2)$.
\end{remark}

\begin{remark}
 $\mathcal O$ may be small but not spherical. For example, for $\g=su(6,\ 3)$ there is a unique orbit $\mathcal O$
 corresponding to the partition $3+3+3$ of $9$. This orbit is small. However, the dimension of a Borel subgroup of
 $\Kc$ is 26, while the dimension of $\mathcal O$ is 27, so that $\mathcal O$ is not spherical.
  (See \cite{colmcg} for an explanation of the parametrization of the
 nilpotent orbits of $su(p,\ q)$.)
\end{remark}
\begin{remark}If $\mathcal O$ is spherical and small, we may have $\kc(x;4)\ne 0$. The principal nilpotent orbit
of $su(2,\ 1)$ satisfies these three conditions.
\end{remark}

\begin{remark}\label{ht2norm}
Suppose that $\mathcal O$ has height 2. Then $\mathcal O$ is
spherical and small.
\end{remark}
\par
\begin{remark}Suppose $\mathcal O$ is spherical and small, then
$\widetilde{V}$ (see Definition \ref{desingdef}) is a commutative
subspace of $\pc$ if and only if $\mathcal O$ has height $2$.
\end{remark}
\medskip
    Let $T$ be a maximal torus in $K$ whose Lie algebra contains $ix$. Let us assume that we have chosen a positive system $\Delta_k^+$ for $(\kc,\ \tc)$ so that $x$ is dominant with respect to $\Delta_k^+$.  $\nk$ is the nilradical of $\kc$ that corresponds to $\Delta_k^+$ and ${\nk}^-$ is the opposite nilradical. Set $u({\mathfrak l}_k) =\lk\cap\nk$. Then, we may assume that $u({\mathfrak l}_k)$ is the
     nilradical of $\lk$, and $u({\mathfrak l}_k)^-$ denotes the opposite nilradical of $\lk$.
\begin{definition} If $\mathcal O$ is a $\Kc$-nilpotent orbit in $\pc$,
and $\lambda\in\widehat{K}$, the set of equivalence classes of irreducible representations of $K$, let $m_{\lambda}(\mathcal O)$ (resp., $m_{\lambda}(\overline{\mathcal O})$) denote the multiplicity of $\lambda$ in $R[\mathcal O]$ (resp., $R[\overline{\mathcal O}]$).
\end{definition}
\begin{definition} We say that $R[\overline{\mathcal O}]$ is self dual as a $K$-module if for every $K$-type $\lambda\in \widetilde{K}$, $m_{\lambda}(\overline{\mathcal O})= m_{{\lambda}^*}(\overline{\mathcal O})$ where ${\lambda}^*\in\widehat{K}$ is the dual of $\lambda$.
\end{definition}
\begin{definition} (Weyl involution of $\gc$) Let $\mathfrak h$ be a Cartan subalgebra of $\gc$ containing $\tc$. Let $\Delta=\Delta(\gc,\ \mathfrak h)$, be the corresponding set of roots. $\Delta^+$ is a positive system for $\Delta$. Assume that we have a Chevalley basis of $\gc$: $\{H_{\alpha},\ X_{\alpha},\ X_{-\alpha}|\ \alpha\in\Delta^+\}$. Then, $\nu$, the Weyl involution of $\gc$ is defined by the requirement that $\nu(H_{\alpha}) = -H_{\alpha}$, and $\nu(X_{\alpha}) = -X_{-\alpha}$ for all $\alpha\in \Delta^+$.
\end{definition}
Let $w_0$ be an element in the normalizer of $T$ in $K$ that represents the longest element of the Weyl group of $(K,\ T)$.
\begin{proposition}\label{pc2red} (Hesselink, Panyushev) As above let $\mathcal O = \Kc\cdot e$ be spherical and let the corresponding triple be $\{x,\ e,\ f\}$. Assume that $\mathcal O$ is small.
Then (1) $\overline{\mathcal O}$ is normal; (2)
$R[\overline{\mathcal O}]^{\nk}$ is a polynomial algebra; (3) If
$R[\overline{\mathcal O}]$ is self dual then there is a
$T$-equivariant ring isomorphism between $R[\overline{\mathcal
O}]^{\nk}$ and $S[\pc(x;2)]^{u({\mathfrak l}_k)}$. (4) If
$R[\overline{\mathcal O}]$ is not self dual, a set of generators
of $R[\overline{\mathcal O}]^{\nk}$ can be identified with
elements in $S(w_0\cdot(\nu(\pc(x;2)))$.
\end{proposition}
\begin{proof} This depends on ideas of Hesselink (from Kempf) applied to the desingularisation of $\overline{\mathcal O}$ in Definition \ref{desingdef} above. Since $\pc(x;i)=0$ if $i>2$, $\widetilde{V}=\pc(x;2)$, and equation(\ref{desingmathcalo}) simplifies to:
\begin{equation}\label{desingsmall}
\pi :\Kc\times _{\Qc\cap \Kc}{\pc(x;2)} \rightarrow
\overline{\mathcal O}=\Kc\cdot\pc(x;2).
\end{equation}
Since $\pc(x;i)=0$ if $i>2$, the nilradical of ${\Qc\cap \Kc}$
acts trivially on $\pc(x;2)$. If we decompose $\pc(x:2)$  into
irreducible $L\cap\Kc$ modules then each of these modules will be
irreducible under $\Qc\cap\Kc$. Thus $\pc(x;2)$ is  completely
reducible as a ${\Qc\cap \Kc}$ module.  From Proposition
\ref{rkctimes}, we know that $R[\Kc\times _{\Qc\cap
\Kc}{\pc(x;2)}]$ is the ring of functions of the normalization and
that $R[\overline{\mathcal O}]\subset R[\Kc\times _{\Qc\cap
\Kc}{\pc(x;2)}]$. Lemma 5.5 of \cite{pan} (based on ideas of
Kempf) implies that:
\begin{equation}\label{normalcond}
R[\Kc\cdot\pc(x;2)]^{{\nk}^-}=R[\Kc\times _{\Qc\cap
\Kc}{\pc(x;2)}]^{{\nk}^-}
\end{equation}
and there is a $T$ equivariant ring homomorphism
\begin{equation}\label{lowestwgt}
 R[\overline{\mathcal O}]^{(\nk)^-}\rightarrow R[\pc(x;2)]^{{u({\mathfrak l}_k)}^-}.
\end{equation}
 Equation (\ref{normalcond}) implies that $R[\overline{\mathcal O}]=R[\Kc\times _{\Qc\cap
\Kc}{\pc(x;2)}]$ and hence $\overline{\mathcal O}$ is normal.\par
$R[\pc(x;2)]^{{u({\mathfrak l}_k)}^-}$ is a polynomial algebra
because $\pc(x;2)$ is a spherical $\lk$-module. (See Corollary
12.2.5 in \cite{gw}.) \par
    $R[\pc(x;2)]$ can be identified with the symmetric algebra $S[\pc(x;2)^*]$,
    where $\pc(x;2)^*$ is the dual of $\pc(x;2)$. But $\pc(x;2)^*$ can be identified with the space $\pc(x;-2)$. The Weyl involution $\nu$
    maps $\pc(x;-2)$ to $\pc(x;2)$. Moreover, this linear isomorphism extends to a ring isomorphism of $S(\pc(x;-2)$ with $S(\pc(x;2)$. This ring isomorphism restricts to a ring isomorphism:
\begin{equation}\label{symringiso}
 S(\pc(x;-2)^{{u({\mathfrak l}_k)}^-}\rightarrow S(\pc(x;2)^{u({\mathfrak l}_k)}.
\end{equation}
The preceding isomorphism carries a lowest weight vector of weight
$\lambda$ to a highest weight vector of weight $-\lambda$.\par The
element $w_0$ defines a ring isomorphism $R[\overline{\mathcal
O}]^{(\nk)^-}\rightarrow R[\overline{\mathcal O}]^{\nk}$. This
fact together with the isomorphisms (\ref{lowestwgt}) and
(\ref{symringiso}) imply that $R[\overline{\mathcal O}]^{\nk}$ is
a polynomial ring.\par
    Now assume $R[\overline{\mathcal O}]$ is self dual, then by normalizing the highest weight and lowest weight vectors in each irreducible $K$ submodule, we can define a ring isomorphism:
\begin{equation}\label{selfdualcase}
 R[\overline{\mathcal O}]^{{\nk}^-}\rightarrow R[\overline{\mathcal O}]^{\nk}.
\end{equation}
The isomorphism takes a lowest weight vector of weight $w_0\lambda$ to a highest weight vector of weight $-w_0\lambda$ in the dual representation.
Combining isomorphisms (\ref{lowestwgt}), (\ref{symringiso}) and (\ref{selfdualcase}), we obtain a $T$-equivariant ring isomorphism
\begin{equation*}
 R[\overline{\mathcal O}]^{\nk}\rightarrow S[\pc(x;2)]^{u({\mathfrak l}_k)}.
\end{equation*}
and thus generators for $R[\overline{\mathcal O}]^{\nk}$ lie in
$S[\pc(x;2)]^{u({\mathfrak l}_k)}$.\par

    When $R[\overline{\mathcal O}]$ is not self dual, it is clear that we can find generators of $R[\overline{\mathcal O}]^{\nk}$ in $S(w_0\cdot(\pc(x;-2)))$.
\end{proof}

    Assume that $\mathcal O$ is spherical and small. Since $R[\pc(x;2)]$ is a unique factorization domain, we can find a set of generators of $R[\pc(x;2)]^{u({\mathfrak l}_k)}$ which are homogeneous and irreducible polynomials.
\begin{definition}\label{genreg} Let $\mathcal O$ be spherical and small.  Set $r$ equal to the $\Kc$-rank of $\mathcal O$. Then $r$ is also the Krull dimension of $R[\overline{\mathcal O}]^{\nk}$. Since $R[\overline{\mathcal O}]^{\nk}$ is a polynomial ring, by Proposition \ref{pc2red}, we can find homogeneous generators $f_{1,\ \mathcal O},\ldots,\ f_{r,\ \mathcal O}$ of $R[\overline{\mathcal O}]^{\nk}$.
Furthermore, if $R[\overline{\mathcal O}]$ is self dual, the
generators can be assumed to be irreducible elements of
$S[\pc(x;2)]^{u({\mathfrak l}_k)}$. Let the degree of $f_{i,\
\mathcal O}$ be denoted by $d_i$ and its $\t$-weight be denoted by
$\gamma_i$.
\end{definition}
\par
\medskip
     The author is not aware of any published work which lists, for each simple $\g$, all the nilpotent $\Kc$-orbits $\mathcal O$ in $\pc$ such that
     $R[\overline{\mathcal O}]$ is self dual. Nor have the spherical orbits been determined for which $R[\overline{\mathcal O}]$ is self
     dual. If $\k$ is semisimple, and all irreducible
     representations of $\k$ are self-dual, then clearly each $R[\overline{\mathcal O}]$ is self
     dual. Thus, $R[\overline{\mathcal O}]$ is self dual for each
     nilpotent $\Kc$-orbit of each of the following simple Lie
     algebras: $su^*(2n)$, $sp(p,\ q)$, $EI$, $EIV$, $EVI$,
      $EVIII$, $EIX$, $FI$, $FII$, and the real split form of
     $G_2$.\par
\begin{remark} For $\g=sl(n,\ \mathbf R)$,  $R[\overline{\mathcal O}]$ is self dual for all spherical
orbits $\mathcal O$ except when $n = 2m$, and $m$ is odd. In that
case, $R[\overline{\mathcal O}]$ is not self-dual if $\mathcal O$
is either of the orbits corresponding to the partition $2^m$.
($2^m$ is the partition of $n$ in which each part size equals
$2$.)  But $R[\overline{\mathcal O}]$ is self dual for all other
spherical orbits.
\end{remark}

\section{Applications to Harish Chandra modules}
    When $\mathcal O$ is spherical and small, Proposition
    \ref{pc2red}, implies that $R[\overline{\mathcal O}]^{\nk}$ is
    a polynomial algebra. Hence the $K$-types $\lambda$ such that $m_{\lambda}(\overline{\mathcal O})\ne
    0$ form a lattice in ${\tc}^*=: Hom(\tc,\ \mathbf C)$. \par
    \begin{definition}
$\Gamma(\overline{\mathcal
O})=\{\lambda\in\widehat{K}|m_{\lambda}(\overline{\mathcal O})\ne
0\}$
    \end{definition}
    In this section we show that if $X$ is an irreducible
    $(\gc,\ K)$-module and $Ass(X)$, the associated variety of $X$, is the
    closure of a nilpotent orbit $\mathcal O$ which is small and
    spherical, then $\Gamma(\overline{\mathcal O})$ determines the
    asymptotic directions of the $K$-types of $X$ in a very
    precise way. (For the definition of $Ass(X)$ see \cite{vo}.)\par
    We first recall some results of Gyoja and Yamashita in \cite{gyya}. Let $\mathfrak w$ be a subalgebra of $\gc$, and ${\mathfrak w}^*$ and ${\gc}^*$ denote their complex dual spaces, i.e., $Hom(\mathfrak w,\ \mathbf C)$ and
$Hom(\gc,\ \mathbf C)$. Let $p:{\gc}^*\rightarrow {\mathfrak w}^*$ denote the projection mapping.
\begin{definition}\label{gycond} (GY-condition) An element $\lambda\in{\gc}^*$ satisfies condition $P_{\kc,\ \mathfrak w}$, if $p(\kc\cdot \lambda) = {\mathfrak w}^*$. Here, $\kc\cdot \lambda=\{z\cdot \lambda| z\in\kc\}$ and ``$\cdot$'' denotes the coadjoint action.
\end{definition}
A subalgebra $\mathfrak w$ for which there is a $\lambda\in{\gc}^*$ satisfying condition $P_{\kc,\ \mathfrak w}$ in Definition \ref{gycond} is  said to satisfy the Gyoja-Yamashita condition (abbreviated GY-condition). Gyoja and Yamashita establish the following result.
\begin{theorem}\label{gythm} (Theorem 2.1 in \cite{gyya}) Let $X$ be an irreducible $(\gc,\ K)$-module, and $\mathfrak w$ be a subalgebra of $\gc$. Suppose that there is an element $\lambda\in Ass(X)$ which satisfies condition $P_{\kc,\ \mathfrak w}$ in Definition \ref{gycond}. Then the action of the enveloping algebra $U(\mathfrak w)$ on $X$ is locally free,i.e., $X$ is a torsion free $U(\mathfrak w)$-module.
\end{theorem}

\begin{definition}\label{formsongc}Let $B_{\gc}$ denote the Killing form of $\gc$. Let $\tau$ denote the conjugation of $\gc$ with respect to the compact real form $\widetilde{u}=\k+i\mathfrak p$. Let $\widetilde{U}$ be the connected subgroup of $\Gc$ with Lie algebra $\widetilde{u}$. Then,
\begin{equation}\label{deflangle}
 \langle z,\ w\rangle :=-B_{\gc}(z,\ \tau(w))
\end{equation}
is a $U$-invariant positive definite Hermitian form on $\gc$.
\end{definition}
    One can verify that $\langle\kc(x;j),\ \pc(x;n)\rangle =0$, for all integers $j$ and $n$, and $\langle(\pc(x;j),\ \pc(x:n)\rangle$ vanishes identically if and only if $j\ne n$. The form $\langle\cdot,\ \cdot\rangle$ in equation (\ref{deflangle}) allows us to identify $\pc(x;i)$ with ${\pc(x;i)}^*$ for each $i$.
\begin{proposition} Suppose that $X$ is an irreducible $(\gc,\ K)$-module,
and $Ass(X)=\overline{\mathcal O}$, where $\mathcal O$ is a
nilpotent $\Kc$-orbit in $\pc$. Assume that $\mathcal O$
corresponds to the normal triple $\{x,\ e,\ f\}$ and that
$\widetilde{V}=\widetilde{V}(\mathcal O):=\pc(x;\ge 2)$ is a
commutative subspace, hence a subalgebra of $\gc$. Then,
$U(\widetilde{V})$ acts locally injectively on $X$.
\end{proposition}
\begin{proof}
Adopt the notation of Definition \ref{desingdef}. By the
observation in Remark \ref{gyyaapp}, $[\qc\cap\kc,\ e] =
\widetilde{V}$. Therefore, the element $e$ satisfies condition
$P_{\kc,\ \widetilde{V}}$ in Definition \ref{gycond}. Thus we can
apply Theorem \ref{gythm} to obtain the desired conclusion.
\end{proof}
\begin{corollary}\label{locinjpc2} Suppose that $X$ is an irreducible $(\gc,\ K)$-module, and $Ass(X)=\overline{\mathcal O}$, where $\mathcal O$ is a nilpotent $\Kc$-orbit in $\pc$ that has height equal to 2. Then, $U(\pc(x;2))$ acts locally injectively on $X$.
\end{corollary}
    Fix a nilpotent orbit $\mathcal O$ which has height equal to 2. Then $\mathcal O$ is spherical. For convenience, assume that $R[\overline{\mathcal O}]$ is self dual. As in Definition \ref{genreg}, let $f_i=f_{i,\ \mathcal O}$, $i=1,\ldots,s$=rank($\mathcal O$) be the generators of  $R[\overline{\mathcal O}]^{\nk}$ which we identify with $S[\pc(x;2)]^{u({\mathfrak l}_k)}$. Let the $\mathfrak t$ weight of $f_i$ be $\gamma_i$. Let $\sigma(f_i)$ be the image of $f_i$ in $U(\gc)$ under the symmetrization map. For each $i$, $\sigma(f_i)\in U(\pc(x;2))$.
    (If $R[\overline{\mathcal O}]$ is not self dual, then we can assume $\sigma(f_i)\in U(w_o\cdot(\pc(x;-2)))$, which also acts locally injectively on
    $X$.)
\begin{corollary} \label{gammamathcalo} Assume the hypotheses of Corollary \ref{locinjpc2} and assume that $R[\overline{\mathcal O}]$ is self dual. Let $\mu$ be a $K$-type in $X$ and let $X^{\mu}$ denote the corresponding highest weight space, then for each $i$, multiplication by $\sigma(f_i)$ maps $X^{\mu}$ to $X^{\mu+\gamma_i}$. Thus, for any choice of non-negative integers $m_1,\ldots, m_s$, $\mu+\sum_{i=1}^s m_i\gamma_i$ is a $K$-type of $X$.
\end{corollary}
In \cite{barvo} Barbasch and Vogan define a notion of asymptotic $K$-support for an irreducible unitary representation of $G$. Let us extend their definition to completely reducible $K$ modules whose $K$-multiplicities are suitably bounded. Assume that $\mathcal R$ is a completely reducible $K$-module with finite $K$-type multiplicities. Let $\mathcal R(\lambda)$ denote the $\lambda$-isotypic component for $\lambda\in\widehat{K}$. If $m(\lambda)=m_{\mathcal R}(\lambda)$ is the multiplicity of $\lambda$ in $\mathcal R$, we will assume that there is some constant $\kappa$ (independent of $\lambda$) such that $m(\lambda)\le \kappa\dim V_{\lambda}$ for all $\lambda\in\widehat{K}$.
\begin{definition}\label{assksupp}(See Theorem 3.6 in \cite{barvo}) \begin{equation*}
\begin{split}
&AS^{K}(\mathcal R)=\{\delta\in{i\t^*}|\exists\ \text{a sequence}\
\delta_n\in\widehat{K},\ n\in{\mathbf N},\ \text{with}\ {\mathcal
R}(\delta_n)\ne 0,\forall n,\cr &\text{and a sequence of positive
numbers},\ t_n,\ \text{with}\lim_{n\rightarrow\infty}t_n=0,\cr
&\text{such that} \lim_{n\rightarrow\infty}t_n\delta_n =\delta\}
\end{split}
\end{equation*}
\end{definition}
We call $AS^{K}(\mathcal R)$ the set of asymptotic directions of
the $K$-types of $\mathcal R$. If $\mathcal O$ is spherical and
small, then ${AS}^K(R[\overline{\mathcal O}])= {\mathbf
R}^+(\Gamma(\overline{\mathcal O}))$, where ${\mathbf
R}^+(\Gamma(\overline{\mathcal O}))$ is the cone spanned by
$\Gamma(\overline{\mathcal O})$ over the nonnegative real numbers.
\par

\begin{theorem}
Suppose that $X$ is an irreducible $(\gc,\ K)$-module, and
$Ass(X)=\overline{\mathcal O}$, where $\mathcal O$ is a nilpotent
$\Kc$-orbit in $\pc$ that is spherical and small. Then,
$AS^{K}(X)={\mathbf R}^+(\Gamma(\overline{\mathcal O}))$.
\end{theorem}
\begin{proof} For convenience of proof, we first assume that $\mathcal O$ has height equal to two,
and $R[\overline{\mathcal O}]$ is self dual. Later we will comment
on removing these restrictions on $\mathcal O$.\par
 To show that ${\mathbf R}^+(\Gamma(\overline{\mathcal O}))\subset AS^{K}(X)$, it suffices to show that $\Gamma(\overline{\mathcal O})\subset AS^{K}(X)$.
 If $\delta\in \Gamma(\overline{\mathcal O})$, and $\mu$ is a $K$-type of $X$, then Corollary \ref{gammamathcalo}
implies that $\mu_n=\mu+n\delta$ is a $K$-type of $X$ for all
non-negative integers $n$. Let $t_n=1/n$. Since
$\lim_{n\rightarrow\infty} t_n\mu_n=\delta$, we have $\delta\in
AS^{K}(X)$.\par
   To show that $AS^{K}(X)\subset {\mathbf R}^+(\Gamma(\overline{\mathcal O}))$, we first choose a good filtration of $U(\gc)$,
   and construct the $(R[\overline{\mathcal O}],\ K)$-module $gr(X)$. It is well known that $gr(X)$ and $X$ are isomorphic as $K$ modules, and
   that $gr(X)$ is a finitely generated
$(R[\overline{\mathcal O}],\ K)$-module. Assume that $y_1,\ldots,\
y_j$ are a set of generators. We can assume that each $y_i$ is a
weight vector for $\mathfrak t$, with weight $\mu_i$. Now assume
that $\delta\in AS^{K}(gr(X))$. Choose corresponding $\t_n$, and
$\delta_n$ as in Definition \ref{assksupp}. Then each $\delta_n$
has the form $\mu_{i_n}+{\delta_n}'$ where ${\delta_n}'$ is a
$K$-type of $R[\overline{\mathcal O}]$. But $||\mu_{i_n}||$  is
bounded as $n\rightarrow \infty$, since the $\mu_{i_n}$ come from
a finite set. It follows that $\delta=\lim_{n\rightarrow\infty}
t_n\delta_n = \lim_{n\rightarrow\infty} t_n{\delta_n}'$. This
implies that $\delta\in {\mathbf R}^+(\Gamma(\overline{\mathcal
O}))$.
\end{proof}

\par
\medskip
\begin{remark} Suppose the $\mathcal O$ is small and spherical, but not of height
two. Then, if $f_i$ is as in Definition \ref{genreg}, one can
still show that the elements $\sigma(f_i)$ act injectively on $X$.
So we still have the analog of Corollary \ref{gammamathcalo}.
\end{remark}

\begin{example} Set $\g=sl(4,\
\mathbf R)$ and $\k=so(4)$. We consider the $K$-type decomposition
of the Speh representation of $SL(4,\ \mathbf R)$, with
Knapp-Vogan parameter $m=-2$, and lowest $K$ type $(1,\ 1)$. We
denote this representation by $S[(-2,\ (1,1)]$, and refer the
reader to pages 586-588 of \cite{knapp} for further discussion of
its properties. Define the following complex symmetric matrices:
\medskip
\begin{equation*}
Y_1=\left
[\begin{array}{rrrr}1&-i&0&0\\-i&-1&0&0\\0&0&0&0\\0&0&0&0\end{array}\right],
Y_2=\left
[\begin{array}{rrrr}0&0&0&0\\0&0&0&0\\0&0&1&-i\\0&0&-i&-1\end{array}\right],
Y_3=\left
[\begin{array}{rrrr}0&0&1&-i\\0&0&-i&-1\\1&-i&0&0\\-i&-1&0&0\end{array}\right].
\end{equation*}
\medskip
 Let $\mathcal O$ be the $\Kc$ orbit of $e=Y_1+Y_2$. $\mathcal O$
is a spherical nilpotent of height 2.\par
     The $K$-decomposition of
$S[(-2,\ (1,1)]$  is as follows:
\begin{equation*}
 S[(-2,\ (1,1)] = \sum_{n=0}^{\infty}\sum_{m=0}^{\infty} \ V(2m+1,\ 2n+1),
\end{equation*}
where $V(2m+1,\ 2n+1)$ is the irr. rep. of $\kc$ with highest
weight $(2m+1,\ 2n+1)$. \par
    The associated variety of $S[(-2,\ (1,1)]$ is
    $\overline{\mathcal O}$. Each of the $Y_i$ ($i=1,2,3$) defined above
    can be viewed as a function on $\pc$, via $Y_i(Z) = trace(Y_iZ)$.
With this identification, the following functions generate
$R[\overline{\mathcal O}]^{\nk}$: $Y_1$ (weight $(2,0)$) and
$Y_1Y_2-(1/4) Y_3^2$ (weight $(2,2)$). The corresponding elements
of $U(\gc)$ act injectively on $S[-2,\ (1,1)]$. The highest
weights of the $K$-types of $S[-2,\ (1,1)]$ may be obtained by
successive application of these elements to the highest weight of
$(1,\ 1)$, the lowest $K$-type of $S[-2,\ (1,1)]$.

\end{example}

\end{document}